\documentclass[12pt
]{article}
\usepackage[T2A]{fontenc}
\usepackage[cp1251]{inputenc}
\usepackage[tbtags]{amsmath}
\usepackage{amsfonts,amsmath,amssymb,mathrsfs,amscd}
\usepackage{amsthm}
\usepackage{verbatim}
\usepackage{graphicx}
\usepackage{xcolor}
\usepackage{enumerate}
\usepackage{textcomp}
\usepackage{latexsym}
\usepackage{verbatim}
\usepackage[mathscr]{eucal}

\theoremstyle{plain} 

\newtheorem*{gluingtheoremA}{Gluing Theorem A}
\newtheorem{gluingtheorem}{Gluing Theorem} 
\newtheorem{proposition}{Proposition}

\theoremstyle{definition}
\newtheorem{definition}{Definition}
\newtheorem{remark}{Remark}

\usepackage[pdftex,unicode,colorlinks,linkcolor=blue,citecolor=red,bookmarksopen,
pdfhighlight=/N]{hyperref}

\newcommand{\RR}{\mathbb{R}} 

\newcommand{\NN}{\mathbb{N}} 
\newcommand{\ZZ}{\mathbb{Z}}

\newcommand{\BB}{\mathbb{B}}  
 
\newcommand{\const}{{\rm const}}

\newcommand{\dd}{\,{\rm d}}

\DeclareMathOperator{\dist}{dist} 
 
\DeclareMathOperator{\clos}{clos} 
\DeclareMathOperator{\Int}{int}
\DeclareMathOperator{\Meas}{Meas}

\DeclareMathOperator{\har}{har}

\DeclareMathOperator{\Har}{har}

\DeclareMathOperator{\sbh}{sbh}

\DeclareMathOperator{\sgn}{sgn}

\DeclareMathOperator{\dom}{dom}
\DeclareMathOperator{\Dom}{Dom}

\title{Gluing Theorems for Subharmonic Functions}
\author{\bf Bulat N. Khabibullin\footnote{This study was financially  supported by the Russian Science Foundation (projects No.~18-01-00002.)}
 \and  \bf Enzhe Menshikova}

\begin{document}

\maketitle
\begin{abstract}
In our articles of recent years, the technique of gluing two subharmonic functions turned out to be very useful in studying the distribution of the roots or masses of holomorphic or subharmonic functions, respectively. Here we develop and improve this technique. Its applications will be given in our further works.



\end{abstract}

\section{Introduction}

As usual, $\mathbb N:=\{1,2, \dots\}$,  $\mathbb R$ and $\mathbb C$ are the sets 
of all {\it natural, real \/} and {\it complex\/} numbers, respectively; 
$\NN_0:=\{0\}\cup\NN$   is French natural series, and $\ZZ:=\NN_0\cup\NN_0$. 

For $d \in \NN$ we  denote by $\mathbb R^d$ the {\it $d$-dimensional real  Euclidean  space\/} with the standard {\it Euclidean norm\/} $|x|:=\sqrt{x_1^2+\dots+x_d^2}$ for $x=(x_1,\dots ,x_d)\in \RR^d$ and the distance function $\dist (\cdot, \cdot)$.

For a subset $S\subset \RR^d$,
 we denote  by $\har (S)$ and  $\sbh (S)$  the classes  of all {\it harmonic\/} (affine for m = 1) and    {\it subharmonic\/} (locally convex for $m = 1$) functions on an open set
$O\supset S$, respectively.  The basis of our note is

\begin{gluingtheoremA}[{\rm \cite[Theorem 2.4.5]{R}, \cite[Corollary 2.4.4]{Klimek}}]\label{gl:th1}
Let $\mathcal O$ be an open set in $\RR^d$, 
and let $\mathcal O_0$ be a  subset of ${\mathcal O}$. If $u\in \sbh({\mathcal O})$, $u_0\in \sbh ({\mathcal O}_0)$, 
and 
\begin{equation}\label{Uu}
\limsup_{y\to x}u_0(y)=u(x) \quad\text{for each $x\in {\mathcal O}\cap \partial {\mathcal O}_0$}, 
\end{equation}
then the formula 
\begin{equation}\label{gU}
U:=\begin{cases}
\max\{u,u_0\} &\text{ on ${\mathcal O}_0$},\\
 u &\text{ on ${\mathcal O}\setminus {\mathcal O}_0$}
\end{cases}
\end{equation}
defines a subharmonic function on ${\mathcal O}$.
\end{gluingtheoremA}
Important applications of Theorem A can be found in our articles \cite{KhaKha19}, \cite{KhaRoz18}, \cite{MenKha19}.

\section{Basic definitions, notations and conventions}\label{Ss12}
\setcounter{equation}{0}

 The reader can skip this Section \ref{Ss12}
and return to it only if necessary.
\subsection{\bf Sets, order, topology.}
For  the {\it real line\/} $\RR$ with  {\it Euclidean norm-module\/} $|\cdot |$,   
\begin{subequations}\label{df:R}
\begin{align}
\RR_{-\infty}:=\{-\infty\}\cup \RR,\; 	\RR_{+\infty}&:=\RR\cup 
\{+\infty\}, \; |\pm\infty|:=+\infty; 
\;  \RR_{\pm\infty}:=\RR_{-\infty}\cup \RR_{+\infty}
\tag{\ref{df:R}$_\infty$}\label{df:Rr}\\
\intertext{is {\it extended real line\/} in the end topology with two ends $\pm \infty$, with   the order relation $\leq$ on $\RR$ complemented by the 
inequalities $-\infty \leq x\leq +\infty$ for $x\in \RR_{\pm\infty}$, with the {\it positive real axis}}
\RR^+:= \{x\in \RR\colon x\geq 0\}, \quad x^+&:=\max\{0,x \},\quad  x^-:=(-x)^+,
\quad \text{for $x\in \RR_{\pm\infty}$},
\tag{\ref{df:R}$^+$}\label{df:R+}
\\
S^+:=\{x\geq 0\colon x\in S \}, \quad S_*&:=S\setminus \{0 \} \quad\text{for $S\subset \RR_{\pm\infty}$}, \quad \RR_*^+:=(\RR^+)_*, 
\tag{\ref{df:R}$_*^+$}\label{df:R*}\\
x\cdot (\pm\infty):=\pm\infty=:&(-x)\cdot (\mp\infty) \quad \text{for $x\in \RR_*^+\cup (+\infty)$}, 
\tag{\ref{df:R}$_\pm$}\label{{infty}+}\\
\frac{x}{\pm\infty}:=0\quad&\text{for $x\in  \RR$,\quad  but $0\cdot (\pm\infty):=0$}
\tag{\ref{df:R}$_0$}\label{{infty}0}
\end{align}
\end{subequations} 
unless otherwise specified. An open connected (sub-)set of $\RR_{\pm\infty}$  is a {\it  (sub-)interval}  of $\RR_{\pm\infty}$. The {\it  Alexandroff\/} one-point {\it compactification\/} of $\mathbb R^d$ is denoted by $\mathbb R^d_{\infty}:=\mathbb R^d \cup \{\infty\}$.

The same symbol $0$ is used, depending on the context, to denote the number zero, the origin, zero vector, zero function, zero measure, etc. 

Given $x\in \RR^d$ and $r\overset{\eqref{df:R+}}{\in} \RR_{+\infty}$, we set 
\begin{subequations}\label{B}
\begin{align}
B(x,r):=\{x'\in \RR^d \colon |x'-x|<r\},&
\quad \overline{B}(x,r):=\{x'\in \RR^d \colon |x'-x|\leq r\},
\tag{\ref{B}B}\label{{B}B}
\\
\quad B(\infty,r):=\{x\in \RR_{\infty}^d \colon |x|>1/r\},&\quad 
\overline B(\infty,r):=\{x\in \RR_{\infty}^d \colon |x|\geq 1/r\},
\tag{\ref{B}$_\infty$}\label{{B}infty}
\\
B(r):=B(0,r),\quad \BB:=B(0,1),& \quad \overline{B}(r):=\overline{B}(0,r),
\quad \overline \BB:=\overline B(0,1).
\tag{\ref{B}$_1$}\label{{B}1}
\\
B_{\circ}(x,r):=B(x,r)\setminus \{x\} ,&\quad 
 \overline{B}_{\circ}(x,r):=\overline{B}(x,r)\setminus \{x\}.
\tag{\ref{B}$_\circ$}\label{Bo}
\end{align}
\end{subequations} 
Thus, the basis of open (respectively closed) neighborhood of the point $x \in \RR_{\infty}^d$ is {\it open\/} (respectively {\it closed\/}) {\it balls\/} 
$B(x,r)$ (respectively $\overline B(x,r)$) centered at  $x$ with radius $r>0$.

Given a subset $S$ of $\RR^d_{\infty}$, the \textit{closure\/} 
$\clos S$, the\textit{ interior\/} $\Int S$  and the \textit{boundary\/} $\partial S$ will always be taken relative $\RR^d_{\infty}$. For $S'\subset S\subset \RR^d_{\infty}$ we write  
$S'\Subset S$ if $\clos S'\subset \Int S$.   
An open connected (sub-)set of $\RR^d_{\infty}$  is a {\it  (sub-)domain}  of $\RR^d_{\infty}$. 

\subsection{\bf Functions.}\label{Functions} Let $X,Y$ are sets. We denote by $Y^X$ the set of all functions  $f\colon X\to Y$. The value $f(x) \in Y$ of an arbitrary function $f\in X^Y$ is not necessarily defined for all $x \in X$. The restriction of a function f to $S \subset X$ is denoted by $f\bigm|_{S}$. We set 
\begin{equation}\label{RX}
\RR_{-\infty}^X\overset{\eqref{df:Rr}}{:=}(\RR_{-\infty})^X, \quad
\RR_{+\infty}^X\overset{\eqref{df:Rr}}{:=}(\RR_{+\infty})^X,\quad
\RR_{\pm\infty}^X\overset{\eqref{df:Rr}}{:=}(\RR_{\pm\infty})^X.  
\end{equation}
A function $f\in \RR_{\pm\infty}^X$ is said to be {\it extended numerical.\/}  
For extended numerical functions $f$, we set 
\begin{equation}\label{dom}
\begin{split}
\Dom_{-\infty}:=f^{-1}(\RR_{-\infty})\subset X,& 
\quad  \Dom_{+\infty} f:=f^{-1}(\RR_{+\infty})\subset X, \\
  \Dom f:=f^{-1}(\RR_{\pm\infty})&=
\Dom_{-\infty} f\bigcup \Dom_{+\infty}f\subset X,\\
\dom f:=f^{-1}(\RR)=&\Dom_{-\infty} f\bigcap \Dom_{+\infty}f\subset X, 
\end{split}
\end{equation} 
For $f,g\in \RR_{\pm\infty}^X$  we write $f= g$ if  
$\Dom f=\Dom g=:D$ and $f(x)=g(x)$ for all $x\in D$, 
and we write $f\leq g$ if $f(x)\leq g(x)$ for all $x\in D$.
For $f\in \RR_{\pm\infty}^X$, $g\in \RR_{\pm\infty}^Y$  and a set $S$, we write 
``$f = g$ {\it on\/} $S$\,'' or  ``$f \leq g$ {\it on\/} $S$\,'' if 
 $f\bigm|_{S\cap D}= g\bigm|_{S\cap D}$ or $f\bigm|_{S\cap D}\leq g\bigm|_{S\cap D}$ respectively.

For $f\in F\subset \RR_{\pm\infty}^X $, we set $f^+\colon x\mapsto \max \{0,f(x)\}$,
$x\in \Dom f$,  $F^+:=\{f\geq 0 \colon f\in F\}$. So, $f$  is \textit{positive\/} on $X$ if $f=f^+$, and  we write ``$f\geq 0$ {\it on\/} $X$''.

The class $\sbh (  O)$  contains the minus-infinity function 
$\boldsymbol{-\infty}\colon x\mapsto -\infty$  identically 
equal to $-\infty$; 
\begin{equation}\label{sbh}
 \sbh_*(  O):=\sbh\,(  O)\setminus 
\{\boldsymbol{-\infty}\}, \quad 
	\sbh^+(  O):=( \sbh (  O))^+.
\end{equation}
 
 If $o\notin O\ni \infty$, then we can to use the  \textit{inversion\/} in  the 
sphere $\partial B(o,1)$ centered at the point $o \in  \RR^d$:
\begin{subequations}\label{stK}
\begin{align}
\star_{o} \colon x\longmapsto x^{\star_{o}}&:= \begin{cases}
o\quad&\text{for $x=\infty$},\\
o+\frac{1}{|x-o|^2}\,(x-o)\quad&\text{for $x\neq o,\infty$},\\
\infty\quad&\text{for $x=o$},
\end{cases}
\qquad \star:=\star_0=:\star_{\infty}
\tag{\ref{stK}$\star$}\label{stK*}
\\
\intertext{together with the  {\it Kelvin transform\/} \cite[Ch. 2, 6; Ch. 9]{Helms}
}
u^{\star_o}(x^{\star_o})&=|x-o|^{d-2}u(x), \quad x ^{\star_o}\in   
O^{\star_o}:=\{x^{\star_o}\colon x\in
  O\}, 
\tag{\ref{stK}u}\label{stKu}
\\
&\Bigl(u\in \sbh (O)\Bigr)\Longleftrightarrow  \Bigl(u^{\star_o}\in \sbh (O^{\star_o})\Bigr).
\tag{\ref{stK}s}\label{stKs}
\end{align}
\end{subequations}

For a subset $S\subset \RR_{\infty}^d$,  the classes $\Har (S)$ and  $\sbh(S)$  consist of the restrictions to $S$ of {\it harmonic\/ {\rm and}  subharmonic functions\/} in some (in general, its own for each function) open set $O\subset \RR_{\infty}^d$ containing $S$.   The class $\sbh_*(S)$ are defined like previous class \eqref{sbh}.

By $\const_{a_1,a_2,\dots}\in \RR$ we denote constants, and constant functions, in general, depend on $a_1,a_2,\dots$ and, unless otherwise specified, only on them,
where the dependence on dimension $d$ of $\RR_{\infty}^d$ will be not specified and not discussed; $\const^+_{\dots}\geq 0$.

\section{General Gluing Theorems}
\setcounter{equation}{0} 

\begin{gluingtheorem}\label{gl:th2}
Let $O$ and $O_0$ be a pair of open subsets in $\RR^d$,  
$v\in \sbh(O)$ and $ v_0\in \sbh(O_0)$ be a pair of  functions such that 
\begin{subequations}\label{g01}
\begin{align}
\limsup_{\stackrel{y\to x}{y\in O_0\cap O}} v(y)&\leq 
v_0(x) \quad\text{for each $x\in O_0\cap \partial O$},
\tag{\ref{g01}$_0$}\label{g010}
\\
\limsup_{\stackrel{y\to x}{y\in O_0\cap O}} v_0(y)&\leq v(x) \quad\text{for each $x\in O\cap \partial O_0$}.
\tag{\ref{g01}$_1$}\label{g011}
\end{align}
\end{subequations}
Then the function 
\begin{equation}\label{Vv}
V:=\begin{cases}
v_0&\text{ on $O_0\setminus O$},\\
\sup \{v_0,v\}&\text{ on $O_0\cap O$},\\
v&\text{ on $O\setminus O_0$,}
\end{cases}
\end{equation}
is subharmonic on $O_0\cup O$.
\end{gluingtheorem}
\begin{proof} It is enough to apply  Gluing Theorem A 
 twice:
\begin{enumerate}
\item[{\bf [O$_0$]}] to one pair of functions  
\begin{subequations}\label{O0}
\begin{align*}
u&:=v_0\in \sbh(O_0), \quad {\mathcal O}:=O_0;
\\
u_0&:=v\bigm|_{O\cap O_0}\in \sbh(O\cap O_0) , 
\quad {\mathcal O}_0:=O\cap O_0\subset O_0, 
\end{align*}
\end{subequations} 
under condition \eqref{g010} realizing condition \eqref{Uu};
\item[{\bf [O]}] 
to another pair of functions 
\begin{subequations}\label{OO}
\begin{align*}
u&:=v\in \sbh(O), \quad {\mathcal O}:=O;
\\
u_0&:=v_0\bigm|_{ O_0\cap O}\in \sbh(O_0\cap O) , 
\quad {\mathcal O}_0:= O_0\cap O\subset O, 
\end{align*}
\end{subequations} 
under condition \eqref{g011} realizing condition \eqref{Uu}.
\end{enumerate}
These two glued subharmonic functions coincide at the open intersection
$O\cap O_0$ and give subharmonic function $V$ defined in \eqref{Vv}.
 
\end{proof}

\begin{gluingtheorem}[{\rm quantitative version}]\label{gl:th3}
Let $O$ and $O_0$ be a pair of open 
subsets in $\RR^d$, and 
$v\in \sbh(O)$ and    $g\in \sbh(O_0)$ be a pair of  functions such that 
\begin{subequations}\label{g01v}
\begin{align}
-\infty <m_v\leq & \inf_{x\in O\cap \partial O_0} v(x), 
\tag{\ref{g01v}m}\label{g01vOm}
\\ 
 \sup_{x\in O_0\cap \partial O}\limsup_{\stackrel{y\to x}{y\in O_0\cap O}} v(y)&\leq M_v<+\infty, 
\tag{\ref{g01v}M}\label{g01vOM}
\\
-\infty < \sup_{x\in O\cap \partial O_0}
\limsup_{\stackrel{y\to x}{y\in O\cap O_0}} g(y)\leq m_g&
< M_g\leq  \inf_{x\in O_0\cap \partial O} g(x) <+\infty.
\tag{\ref{g01v}g}\label{g01g}
\end{align}
\end{subequations}
If we choose the function 
\begin{equation}\label{v0g}
v_0:=\frac{M_v^++m_v^-}{M_g-m_g} (2g-M_g-m_g)\in \sbh(O_0), 
\end{equation} 
then the function $V$ from \eqref{Vv} is subharmonic on $O_0\cup O$. 
\end{gluingtheorem}
\begin{proof}
The function $v_0$ from definition \eqref{v0g} is subharmonic on $O_0$ since this function $v_0$  has a form $\const^+g+\const$ with $\const^+\in \RR^+$, $\const \in \RR$. 
In addition, by construction \eqref{v0g}, for each $x\in O_0\cap \partial O$, we obtain
\begin{multline*}
\limsup_{\stackrel{y\to x}{y\in O_0\cap O}} v(y)\overset{\eqref{g01vOm}-\eqref{g01vOM}}{\leq}
 M_v^++m_v^-= \frac{M_v^++m_v^-}{M_g-m_g} (2M_g-M_g-m_g)\\
 \overset{\eqref{g01g}}{\leq}\frac{M_v^++m_v^-}{M_g-m_g} \Bigl(2 \inf_{x\in O_0\cap \partial O} g(x)  -M_g-m_g\Bigr)\\
=\inf_{x\in O_0\cap \partial O} \frac{M_v^++m_v^-}{M_g-m_g} \Bigl(2  g(x)  -M_g-m_g\Bigr)\\
\overset{\eqref{v0g}}{=} \inf_{O_0\cap\partial O} v_0\leq 
v_0(x), \quad \forall x \in O_0\cap \partial O.
\end{multline*} 
Thus,  we have  \eqref{g010}. 
Besides, by construction \eqref{v0g}, for each $x\in O\cap \partial O_0$, we obtain
\begin{multline*}
\limsup_{\stackrel{y\to x}{y\in O_0\cap O}} v_0(y)
\overset{\eqref{v0g}}{\leq}
\frac{M_v^++m_v^-}{M_g-m_g} \biggl(2\limsup_{\stackrel{y\to x}{y\in O_0\cap O}} g-M_g-m_g\biggr)
\overset{\eqref{g01g}}{\leq}
\frac{M_v^++m_v^-}{M_g-m_g} (2m_g-M_g-m_g)
\\=-(M_v^++m_v^-)\leq  -m_v^-\leq m_v
\overset{\eqref{g01vOm}}{\leq}
 \inf_{x\in O\cap \partial O_0} v(x)\leq v(x), 
 \quad \forall x\in  O\cap \partial O_0.
\end{multline*} 
Thus,  we have  \eqref{g011}, and our Gluing Theorem \ref{gl:th3} follows from 
Gluing Theorem \ref{gl:th2}. 
\end{proof}

\begin{remark} 
Theorems of this section can be easily transferred to the cone of plurisubharmonic functions
\cite[Corollary 2.9.15]{Klimek}.
We sought to formulate our theorems and their proofs with the possibility of their fast transport to the plurisubharmonic functions and to abstract potential theories with  more general constructions based on the theories of harmonic spaces and sheaves  in the spirit 
of books  \cite{BB80}, \cite{BB66},  \cite{BBC81}, \cite{BH}, \cite{AL},   etc.
\end{remark}

\section{Gluing with the Green Function}
\setcounter{equation}{0}

\begin{definition}[{\rm \cite{R}, \cite{HK}, \cite{Landkoff}}]\label{df:kK} 
 For $q\in \RR$, we set  
\begin{subequations}\label{kK}
\begin{align}
k_q(t)& := \begin{cases}
\log t  &\text{ if $q=0$},\\
 -\sgn (q)  t^{-q} &\text{ if $q\in \RR_*$,} 
\end{cases}
\qquad  t\in \RR_*^+,
\tag{\ref{kK}k}\label{{kK}k}
\\
K_{d-2}(x,y)&:=\begin{cases}
k_{d-2}\bigl(|x-y|\bigr)  &\text{ if $x\neq y$},\\
 -\infty &\text{ if $x=y$ and $d\geq 2$},\\
0 &\text{ if $x=y$ and  $d=1$},\\
\end{cases}
\quad  (x,y) \in \RR^d\times \RR^d.
\tag{\ref{kK}K}\label{{kK}K}
\end{align}
\end{subequations}
 \end{definition}

Reminder, that a set $E\subset \RR^d$ is called {\it polar\/} if there is a function $u\in \sbh_*(\RR^d)$ such that 
\begin{equation}\label{E}
\Bigl(E\subset (-\infty)_u:=\{ x\in \RR^d\colon u(x)=-\infty \} \Bigr)
\Longleftrightarrow \Bigl(\text{Cap}^* E=0\Bigr),
\end{equation}
where the set $(-\infty)_u$ is {\it minus-infinity\/} $G_{\delta}$-set for the function 
$u$,   
\begin{equation*}
\text{Cap}^*(S):=
\inf_{S\subset O=\Int O}  
\sup_{\stackrel{C=\clos C\Subset O}{\mu\in \Meas^{1+}(C)}} 
 k_{d-2}^{-1}\left(\iint K_{d-2} (x,y)\dd \mu (x) \dd \mu(y)\right) 
\end{equation*}
is the {\it outer capacity\/} of $S\subset \RR^d$.

Let $\mathcal O$ be an \textit{open  proper} subset in $\RR_{\infty}^d$. 

Consider a point $o\in \RR^d$ and subsets $S_0, S\subset \RR_{\infty}^d $  such that  
  \begin{equation}\label{x0S}
\RR^d \ni o \in \Int S_0\subset S_0\Subset S \subset \Int \mathcal O
=\mathcal O\subset \RR^d_{\infty}\neq \mathcal O. 
\end{equation} 
Let $D$ be a \textit{domain\/} in $\RR_{\infty}^d$ such that
\begin{equation}\label{Dg}
 o\overset{\eqref{x0S}}{\in} \Int  S_0\subset S_0\Subset D \Subset S \subset \mathcal O.
\end{equation}
Such domain $D$ possesses the generalized Green's function $g_D (\cdot, o)$\/  {\rm (see \cite[5.7.2]{HK}, \cite[Ch.~5, 2]{Helms})} with pole at the point 
$o\overset{\eqref{Dg}}{\in} D$ 
 described by the following properties:  
\begin{subequations}\label{gD}
\begin{align}
g_D(\cdot , o)&\in \sbh^+ \bigl(\RR_{\infty}^d\setminus \{o\}\bigr)\subset \sbh^+\bigl(\mathcal O\setminus \{o\}\bigr) , 
\tag{\ref{gD}s}\label{gDs}\\
g_D(\cdot ,o)&= 0\text{ on $\RR_{\infty}^d\setminus \clos D \supset \mathcal O\setminus \clos D\supset \mathcal O\setminus S$}, 
\tag{\ref{gD}$_0$}\label{gD0}\\
g_D(\cdot , o)&\in \har \bigl(D\setminus \{o\}\bigr)\subset 
\har\bigl(S_0\setminus\{o\}\bigr)\overset{\eqref{Bo}}{\subset} \har\bigl(B_{\circ}(o,r_o)\bigr)
\tag{\ref{gD}h}\label{gDh}\\
\intertext{with a number $r_o\in \RR_*^+$, $g_D(o,o):=+\infty$,}
g_D(x,o)&\overset{\eqref{{kK}K}}{=}-K_{d-2}(x,o)+O(1) \quad\text{when $o\neq x\to o$}.
\tag{\ref{gD}o}\label{gD0a}\\
\intertext{and the following strictly positive number}
0<M_g&:=\inf_{x\in \partial S_0} g_D (x,o)=
\const^+_{o, S_0, D,S} 
 \tag{\ref{gD}M}\label{Mg}\\
\intertext{depends only on $S_0,S,D$ and the pole $o$, and, by the minimum principle, we have}
g_D(x,o)-M_g&\overset{\eqref{Mg}}{\geq} 0\quad\text{for all $x\in S_0\setminus \{o\}$}.
\tag{\ref{gD}M+}\label{Mg+}
\end{align}
\end{subequations}
Properties \eqref{gD} for the generalized Green's function 
$g_D (\cdot, o)$\/  
from \eqref{x0S}--\eqref{Dg} are well known \cite[4.4]{R}, \cite[5.7]{HK}, and 
property \eqref{x0S} follows from $0< g(\cdot ,o)\in C\bigl(D\setminus \{o\}\bigr)$ on $D\setminus \{o\}$.  

\begin{gluingtheorem}\label{gl:th4}
Under conditions  \eqref{x0S} 
suppose that a function $v \in \sbh(\mathcal O\setminus S_0)$ satisfy constraints above and below  in the form
\begin{equation}\label{vabS}
-\infty<m_v\overset{\eqref{g01vOm}}{\leq} \inf_{S\setminus S_0} v\leq  
\sup_{S\setminus  S_0} v\overset{\eqref{g01vOM}}{\leq}  M_v<+\infty. 
\end{equation}
Every domain $D$ with inclusions  \eqref{Dg} possesses  the generalized Green's function $g_D(\cdot , o)$ with pole $o\in \Int S_0$, properties  \eqref{gD} and constant $M_g$ of   \eqref{Mg} such that the choice of function 
\begin{subequations}\label{gDV}
\begin{align}
v_0&\overset{\eqref{v0g}}{:=}\frac{M_v^++m_v^-}{M_g} \bigl(2g_D(\cdot , o) -M_g\bigr)\in 
\sbh\bigl(\RR_{\infty}^d\setminus \{o\}\bigr)\subset 
\sbh\bigl(\Int S\setminus \{o\}\bigr),
\tag{\ref{gDV}v}\label{gDVv}
\\
\intertext{defines the subharmonic function}
V&\overset{\eqref{Vv}}{:=}\begin{cases}
v_0&\text{ on $S_0$},\\
\sup \{v_0,v\}&\text{ on $S\setminus S_0$},\\
v&\text{ on $\mathcal O\setminus S$,}
\end{cases} \qquad \text{from $\sbh_*\bigl(\mathcal O\setminus \{o\}\bigr)$,}
\tag{\ref{gDV}V}\label{gDVV}\\
\intertext{satisfying the conditions}
 V&\overset{\eqref{gDh}}{\in} \har \bigl(S_0\setminus \{o\}\bigr)\overset{\eqref{Bo}}{\subset} \har\bigl(B_{\circ}(o,r_o)\bigr)
\quad\text{with a number $r_o\in \RR_*^+$},
\tag{\ref{gDV}h}\label{gDhV}\\
V&\overset{\eqref{Mg}-\eqref{Mg+}}{\geq 0} \quad\text{on $ S_0$,}
\tag{\ref{gDV}+}\label{gDhV+}\\
V(x)&\overset{\eqref{gD0a}}{=}-2\frac{M_v^++m_v^-}{M_g} K_{d-2} (x,o)+O(1) \quad\text{when $o\neq x\to o$}.
\tag{\ref{gDV}o}\label{gD0aV}
\end{align}
\end{subequations}
\end{gluingtheorem}
\begin{proof} It is enough to apply Gluing Theorem \ref{gl:th3}
 with 
\begin{equation*}
O:=\mathcal O\setminus \clos S_0, 
\quad O_0:=\Int S\setminus \{o\}, \quad 
g:=g_D(\cdot,o), \quad  m_g:=0 
\end{equation*} 
according to the references written above the relations in 
 \eqref{vabS}--\eqref{gDV}.
\end{proof}

Given $S\subset \RR^d$ and $r\in \RR^+$, a set 
\begin{equation}\label{Scup}
S^{\cup r}:=S\bigcup \bigcup\limits_{x\in S} B(x,r).
\end{equation}
is called a outer $r$-parallel set \cite[Ch.~I,\S~4]{Santalo}. 
Easy to install the following
\begin{proposition}\label{llemmaD} 
Let a subset $S_0\Subset \RR^d$
 be connected, and $r\in \RR_*^+$. Then $S^{\cup r}$ is connected, $S_0\Subset S_0^{\cup (r/3)}$, 
 and there is a regular for the Dirichlet problem domain $D\subset S_0^{\cup r}$ such 
that
\begin{equation}\label{SDS}
S_0^{\cup (r/3)}\Subset \partial D\Subset S_0^{\cup (2r/3)}. 
\end{equation}
\end{proposition}

For $v\in L^1\bigl(\partial B(x,r)\bigr)$, we define the averaging value of $v$ at the point $x$ as
\begin{equation}\label{v0}
v^{\circ r}(x):=\frac{1}{}\int_{\partial B(x,r)} v\dd \sigma_{d-1},
\end{equation}
 where $\sigma_{d-1}$ the normalized by $1$ surface measure on the sphere $\partial B(x,r)$.

\begin{gluingtheorem}\label{gl:th_es}
Let $\mathcal O\subset \RR^d$ be an open subset, and  $S_0\subset \RR^d$ be a connected set such that there is a point
\begin{equation}\label{S0}
 o\overset{\eqref{x0S}}{\in} \Int  S_0\subset S_0 \Subset \mathcal O.
\end{equation}
Let  $r\in \RR^+$ be a  number such that
\begin{equation}\label{posr}
0<r<\dist(S_0, \partial \mathcal O),
\end{equation}
and $D$ be a domain from Proposition\/ {\rm \ref{llemmaD}} satisfying \eqref{SDS}.
Let $v \in \sbh_*(\mathcal O\setminus S_0)$ be a function satisfying constraints above and below  in the form
\begin{subequations}\label{avv}
\begin{align}
v&\leq M_v <+\infty \quad\text{on $S_0^{\cup r}\setminus S_0$},
\tag{\ref{avv}M}\label{avvM}
\\
m_v&:=\inf \bigl\{ v^{\circ (r/3)}(x)\colon x\in S_0^{\cup (2r/3)}\setminus S_0^{\cup (r/3)}\bigr\}.
\tag{\ref{avv}m}\label{avvm}
\end{align}
\end{subequations}
Then $m_v>-\infty $, and   there is a function $V\in \sbh_*(\mathcal{O}\setminus \{o\}) $ satisfying \eqref{gDhV}--\eqref{gDhV+}, i.\,e., 
\begin{subequations}\label{VK}
\begin{align}
V&\in \har^+\bigl(S_0\setminus \{o\}\bigr),
\tag{\ref{VK}h}\label{gDVh}\\
V&=v \quad \text{on $\mathcal O\setminus \clos S^{\cup r}$},
\tag{\ref{VK}=}\label{gDV=}\\
\intertext{and such that}
V(x)&\overset{\eqref{gD0aV}}{=}- 2\frac{M_v^++m_v^-}{M_g}K_{d-2} (x,o)+O(1) \quad\text{when $o\neq x\to o$}
\tag{\ref{VK}o}\label{gDVo},\\
\intertext{where}
0<M_g&:=\inf_{x\in \partial S_0^{\cup (r/3)}} g_D (x,o)=
\const^+_{o, S_0, r,D}
\tag{\ref{VK}g}\label{gDVg}
\end{align}
\end{subequations}

\end{gluingtheorem}
\begin{proof}
 We have $m_v>-\infty$ since the function  
$v^{\circ (r/3)}$  is continuous in $ S_0^{\cup (2r/3)}\setminus S_0^{\cup (r/3)}$
 \cite[Theorem 1.14]{Helms}. 
The function  $v$ can be transformed using the Perron\,--\,Wiener\,--\,Brelot method (into the open ``layer'' $S_0^{\circ r}\setminus \clos S_0$ from 
boundary of this layer) to a new \textit{subharmonic function $\widetilde v\geq v$ on} $\mathcal O \setminus S_0$ such that $\widetilde v\in \har (S_0^{\circ r}\setminus \clos S_0)$ and 
$\widetilde v=v$ on $\mathcal O\setminus \clos S$.
This follow from the principle of subordination (domination) for harmonic continuations
and  the  maximum principle that 
\begin{equation}\label{estwv}
-\infty <m_v\leq \widetilde v\quad\text{on $S_0^{\cup (2r/3)}\setminus S_0^{\cup (r/3)}$}, \quad
\quad \widetilde v\leq M_v\quad\text{on $S_0^{\cup r}\setminus S_0$},
\end{equation}
If we choose in Gluing  Theorem \ref{gl:th4} for the role a set $S_0$ the set $S_0^{\cup (r/3)}$, and instead of $S$ the set $S_0^{\cup (2r/3)}$, then, by construction \eqref{gDVv}--\eqref{gDVV} and conditions \eqref{gDhV}--\eqref{gD0aV} , we get 
series of conclusions \eqref{VK} of Theorem \ref{gl:th_es}.
\end{proof}

\vskip 1 mm 

\textit{Russian Federation, Bashkortostan, Ufa, Bashkir State University}

E-mail: \textit{Khabib-Bulat@mail.ru, algeom@bsu.bashedu.ru}

\end{document}